\documentclass[12pt]{article}
\usepackage{amsmath,amsthm,graphics}

\theoremstyle{plain}
\newtheorem{theorem}{Theorem}
\newtheorem{lemma}{Lemma}

\newcommand\A{{\mathcal A}}

\newcommand\T{{\mathcal T}}
\newcommand\C{{\mathcal C}}

\title{Self-describing sequences and the Catalan family tree}

\author{Zoran \v Suni\'k\\
\small Department of Mathematics\\[-0.8ex]
\small Texas A\&M University\\[-0.8ex]
\small College Station, TX 77843-3368, USA}

\date{\small MR Subject Classifications: 05A15, 05C05, 11Y55}

\begin{document}
\maketitle

\begin{abstract}
We introduce a transformation of finite integer sequences, show
that every sequence eventually stabilizes under this
transformation and that the number of fixed points is counted by
the Catalan numbers. The sequences that are fixed are precisely
those that describe themselves --- every term $t$ is equal to the
number of previous terms that are smaller than $t$. In addition,
we provide an easy way to enumerate all these self-describing
sequences by organizing them in a Catalan tree with a specific
labelling system.
\end{abstract}

\subsection*{Prefix ordered sequences and rooted labelled trees}
The following connection between prefix ordered sequences and
rooted labelled trees is well known and we briefly mention only
the instance which is useful for our considerations.

\begin{figure}[!ht]
  \begin{center}
  \includegraphics{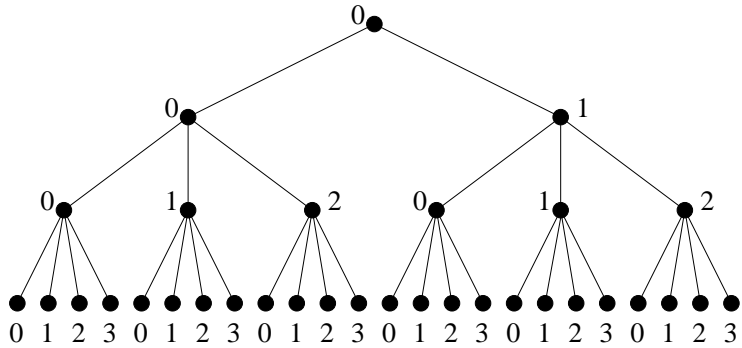}
  \end{center}
  \caption{The rooted labelled tree $\T$ up to the third generation}
  \label{factorialtree}
\end{figure}

Let $\A$ be the set of finite integer sequences
$a=(a_0,a_1,\dots)$ with the property that $0 \leq a_i \leq i$,
for all indices. We order the sequences in $\A$ by the
\emph{prefix} relation, i.e.,
\[ (a_0,a_1,\dots,a_n) \preceq (b_0,b_1,\dots,b_m) \]
if $n \leq m$ and $a_i=b_i$, for $i=0,\dots,n$. The sequences in
$\A$ can be organized in a rooted labelled tree $\T$ which
reflects the prefix order relation. The root of the tree $\T$ is
labelled by $0$. Every vertex that is at distance $n$ from the
root has $n+2$ children labelled by $0,1,\dots,n,n+1$ (see
Figure~\ref{factorialtree}).
The vertices whose distance to the root is $n$ form the $n$-th
\emph{level} of the tree $\T$, which is also called the $n$-th
\emph{generation}. For every vertex $v$ at the level $n$ in the
tree $\T$ there exist a unique path of length $n$ from the root to
$v$. The labels of the vertices on this path form a unique
sequence $(a_0,a_1,\dots,a_n)$ in $\A$ that corresponds to the
vertex $v$ and this sequence is called the \emph{full name} of
$v$. The correspondence
\[ v \leftrightarrow \text{the full name of }v \]
provides a bijection between the vertices in $\T$ and the
sequences in $\A$. Under this bijection, the vertices from the
$n$-th generation in $\T$ correspond to the sequences of length
$n+1$ in $\A$. The set of vertices in the $n$-th generation is
denoted by $\T_n$ and the corresponding set of sequences by
$\A_n$.

The sequence $a=(a_0,a_1,\dots,a_n)$ is a prefix of the sequence
$b=(b_0,b_1,\dots,b_m)$ if and only if the vertex $v_a$ with full
name $a$ is on the unique path between the root and the vertex
$v_b$ with full name $b$, i.e., if and only if the vertex $v_a$ is
an ancestor of the vertex $v_b$. Consider a graph endomorphism
$\alpha$ of $\T$ that fixes the root (and therefore also preserves
the levels). Such an endomorphism corresponds to a transformation
of sequences $\alpha:\A \to \A$ that preserves the length of the
sequences and also their prefix order, i.e.,
\[ a \preceq b \qquad \text{implies} \qquad \alpha a \preceq
    \alpha b, \]
for all sequences $a$ and $b$ in $\A$.

In the sequel, we often deliberately blur the distinction between
the vertices in $\T$ and the corresponding sequences in $\A$.
Similarly, we do not distinguish tree endomorphisms of $\T$ fixing
the root from sequence transformations that preserve the length
and the prefix order. This mistake actually improves our
presentation.

Let $\alpha$ be an endomorphism of $\T$. Since every generation in
$\T$ is finite, the $\alpha$ \emph{orbit}
\[ \alpha^*u = \{ \alpha^iu \;|\; i \geq 0 \;\} \]
of every vertex $u$ of $\T$ is finite. Thus, starting from any
vertex, repeated applications of $\alpha$ produce \emph{periodic
points}, i.e., points $a$ for which $\alpha^k a=a$ for some $k>
0$. The \emph{period} of the periodic point $a$ is the smallest
$k$ for which $\alpha^k a=a$. The points of period 1 are
\emph{fixed points} and the points of period dividing $2$ are
\emph{double points}. Obviously, if $u$ and $v$ are periodic
points of $\alpha$ and $u$ is a prefix of $v$ then the period of
$u$ divides the period of $v$.

It is easy sometimes to estimate how long it takes before a
periodic point is reached. We make use of the
\emph{lexicographic ordering} $\leq$ of the sequences in $\A_n$
(note the difference with the prefix ordering $\preceq$).
Namely, for $a=(a_0,a_1,\dots,a_n)$ and $b=(b_0,b_1,\dots,b_n)$,
set $a<b$ if $a_i<b_i$ at the first index where $a$ and $b$
differ.

\begin{theorem}\label{on2}
Let $\alpha$ be an endomorphism of the tree $\T$ and assume that,
for some $n \geq 1$, there exists $k \geq 1$ such that, for every
vertex $u$ in generation $n$, either
\[ u \leq \alpha^k u \leq \alpha^{2k} u \leq \dots \]
or
\[ u \geq \alpha^k u \geq \alpha^{2k} u \geq \dots .\]
Then, starting from any point in generation $n$, repeated
applications of $\alpha$ lead to a periodic point of period
dividing $k$ in $O(n^2)$ steps.
\end{theorem}
\begin{proof}
We show that $\beta=\alpha^k$ reaches a fixed point in no more
than
\[ 1 + 2 + \dots + n = n(n+1)/2 \]
steps.

Start with any vertex $u$ in generation $n$. Without loss of
generality we may assume
\[ u \leq \beta u  \leq \beta^2 u \leq \dots .\]
After the first application of $\beta$ the initial segment up to
index 1 of $\beta u$ is fixed under $\beta$. After the next two
steps the entry at index 2 will be fixed. Proceeding in the same
fashion we see that the initial segment of $\beta^{1+2+ \dots
+i} u$ up to index $i$ is fixed under $\beta$. Indeed, once the
initial segment up to index $i-1$ is fixed the entry at index
$i$ can go up no more than $i$ times (from $0$ to $i$) before it
stabilizes. Thus, $\beta^{1+2+ \dots +n}u$ is fixed under
$\beta$.
\end{proof}

\subsection*{Self-describing sequences}
We define an endomorphism $\delta: \A \to \A$ transforming
sequences in $\A$ by
\[ (\delta a)_i = \# \{j\;|\; j < i,\; a_j < a_i \}. \]
Thus, for each term $t$ in the sequence $a$, $(\delta a)_i$ counts
the number of previous terms that are smaller than $t$. The
transformation $\delta$ makes perfect sense even for sequences out
of $\A$, but the image is in $\A$ and it stays there under further
iterations. A sequence that is fixed under $\delta$ is called a
\emph{self-describing sequence}. Therefore, the sequence
$a=(a_0,a_1,\dots)$ is self-describing if
\[ \# \{j\;|\; j < i,\; a_j < a_i \} = a_i, \]
for all indices, i.e., every term $t$ is equal to the number of
previous terms that are smaller than $t$.

\subsection*{The Catalan family tree} 
\begin{figure}[!ht]
  \begin{center}
   \includegraphics{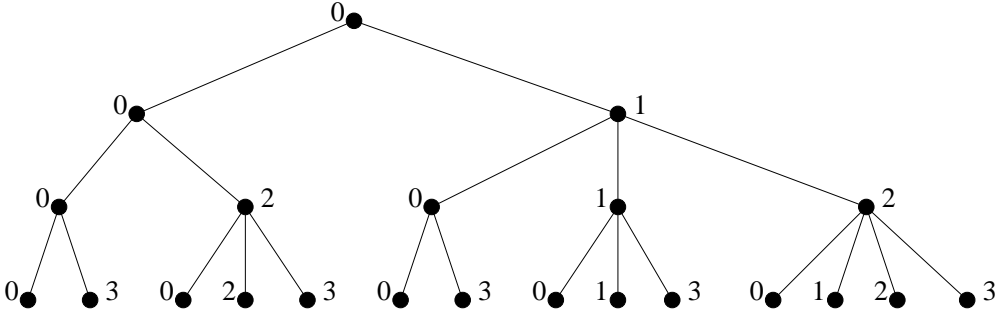}
  \end{center}
  \caption{The Catalan family tree $\C$ up to the third generation}
  \label{ctlsds}
\end{figure}

We describe now a rooted
labelled subtree of $\T$, denoted by $\C$ and called \emph{the
Catalan family tree} or just the \emph{Catalan family}. The root
vertex $0$ belongs to $\C$. It has two children named $0$ and $1$
and we consider $0$ the older sibling. The oldest sibling in this
family always has 2 children, the second oldest 3, the third
oldest 4, and so on. The oldest child of a member of the family
$x$ gets named after the oldest sibling of $x$, the second oldest
child after the second oldest sibling, and so on, until $x$ uses
its own name for its second to last child and $n$ for the youngest
one, where $n$ is the generation number of the children (the level
in the tree). The diagram in Figure~\ref{ctlsds} depicts the
family members of $\C$ up to the third generation.

\subsection*{The connection} We establish now a connection between
the self-describing sequences and the Catalan family tree.

\begin{theorem}
The full names of the members of the Catalan family are precisely
the self-describing sequences. In other words, they are the fixed
points of the endomorphism $\delta$.

Moreover, repeated applications of $\delta$ to any sequence in
$\A$ eventually produce a member of the Catalan family, i.e. a
fixed point of $\delta$. The number of applications needed to
reach such a point is $O(n^2)$.
\end{theorem}

All statements of the theorem are implied by Theorem~\ref{on2} and
the following lemma.

\begin{lemma}
If $a$ is a member of the Catalan family then $a = \delta a$.
Otherwise, $a < \delta a$.
\end{lemma}

\begin{proof}
The proof is by induction on the generation number $n$. The
statement is true for $n=0$ and $n=1$. Assume that the statement
is true for all vertices up to the $n$-th generation.

Let
\[ a=(a_0,a_1,\dots,a_n,x) \]
be a $(n+1)$-st generation member of the Catalan family. We
consider two cases.

If $x=n+1$ then
\[ \# \{j\;|\; j < n+1,\; a_j < x \} =
   \# \{j\;|\; j < n+1,\; a_j < n+1 \} = n+1 = x,  \]
and $a$ is a fixed point of $\delta$.

If $x \neq n+1$, then $a_n \geq x$ and there exists an $n$-th
generation member of the Catalan family whose full name is
\[ a'=(a_0,a_1,\dots,a_{n-1},x), \]
namely the one after whom $a$ was named. We have
\[ \# \{j\;|\; j < n+1,\; a_j < x \} =
   \# \{j\;|\; j < n,\; a_j < x \} = x,\]
where the first equality comes from the fact that $a_n \geq x$ and
the second from the inductive hypothesis, since $\delta a'=a'$.

Thus all members of the Catalan family are fixed under $\delta$.

Now, let
\[ a=(a_0,a_1,\dots,a_n,x) \]
be a full name of a vertex in $\T$ in the $n$-th generation that
is not a member of the Catalan family $\C$. If any proper prefix
of $a$ is not in $\C$ we obtain the claim directly from the
inductive hypothesis. Thus we may assume that
\[ a''=(a_0,a_1,\dots,a_n) \]
is a member of the Catalan family. Since $a$ is not in $\C$ we
have $a_n \neq x$ and $n+1 \neq x$. We consider two cases.

If $a_n > x$ then $a'=(a_0,a_1,\dots,a_{n-1},x)$ is not in $\C$
and
\[ \# \{j\;|\; j < n+1,\; a_j < x \} =
   \# \{j\;|\; j < n,\; a_j < x \} > x,\]
where the equality comes from the fact that $a_n > x$ and the
inequality from the inductive hypothesis.

If $a_n < x < n+1$ then
\[ \# \{j\;|\; j < n+1,\; a_j < x \} =
   \# \{j\;|\; j < n,\; a_j < x \} + 1 \geq x+1,\]
where the equality comes from the fact that $a_n < x$ and the
inequality from the inductive hypothesis. The equality in the last
case is possible only when $a'=(a_0,a_1,\dots,a_{n-1},x)$ is in
$\C$.
\end{proof}

We proceed by counting the self-describing sequences with fixed
length. In addition, we obtain a result on the distribution of
names in $\C$. Recall that the $n$-th Catalan number is equal to
\[c_n = \frac{1}{n+1}\binom{2n}{n}.  \]
A recursive definition of the Catalan numbers is given by
\begin{gather*}
 c_0 = 1,\\
 c_{n+1} = c_0c_n + c_1c_{n-1} + \cdots + c_nc_0.
\end{gather*}

\begin{theorem}\label{thm:count}
The number of self-describing sequences in $\A_n$, i.e., the
number of $n$-th generation members of the Catalan family is the
$(n+1)-th$ Catalan number $c_{n+1}$.

Moreover, for $r=0,\dots,n$, the number of $n$-th generation
members of the Catalan family whose name is $r$ is equal to
$c_rc_{n-r}$.
\end{theorem}
\begin{proof}
Denote by $z_n$ the number of $n$-th generation members of the
Catalan family whose name is $0$. More generally, for
$r=0,\dots,n$ denote by $f_{n,r}$ the number of $n$-th generation
members of the Catalan family whose name is $r$. Finally, denote
by $g_n$ the number of $n$-th generation members of the Catalan
family.

Since the oldest child of every member of the Catalan family is
named $0$, we have, for all $n$,
\[ z_{n+1} = g_n. \]

Since the youngest sibling in the $r$-th generation is always
named $r$ and the oldest $0$ we also have, for all $r$,
\[ f_{r,r} = f_{r,0} = z_r. \]

For some fixed $r$, consider the set of $f_{r,r}$ $r$-th
generation members named $r$ together with all their descendants
in $\C$ whose names are greater or equal to $r$. This forest of
$f_{r,r}$ identical subtrees of $\C$ contains all members of
$\C$ whose name is $r$. Moreover, each tree in this forest looks
exactly like the Catalan family tree, except that all labels are
increased by $r$. Indeed, each $r$-th generation member of $\C$
named $r$ has two children, named $r$ and $r+1$, the oldest
sibling always has two children, the second oldest three, etc.
Thus, for any $n$ and $r=0,\dots,n$, the number $f_{n,r}$ of
$n$-th generation members of $\C$ named $r$ is $f_{r,r}$ times
larger than the number of $(n-r)$-th generation members of $\C$
named $0$, i.e.,
\[ f_{n,r} = f_{r,r} f_{n-r,0} = z_r z_{n-r}. \]

Since $z_0=1$ and
\begin{align*}
 z_{n+1} &= g_n = f_{n,0} + f_{n,1} + \cdots + f_{n,n}\\
         &= z_0z_n + z_1z_{n-1} + \cdots + z_nz_0
\end{align*}
we conclude that, for all $n$, $z_n$ is the $n-th$ Catalan number.
The statements of the theorem follow now easily from the relations
$g_n=z_{n+1}$ and $f_{n,r} = z_r z_{n-r}$.
\end{proof}

\subsection*{Connection to other Catalan trees and objects} It is
well known that the Catalan numbers appear naturally under many
circumstances. The exercises on Catalan numbers in
\cite{stanley:ec2} provide a trove of examples, along with
references, in which Catalan numbers count the number of objects
of particular type and size. The self-describing sequences provide
yet another example that we now relate to some other objects
counted by the Catalan numbers.

Consider the sequences in $\A$ with the property that $a_{i+1}
\leq a_i+1$, for all indices (see the Exercise 6.19.u in
\cite{stanley:ec2}). Such sequences are called \emph{sequences
with unit increase}. The rooted labelled tree that corresponds to
the set of sequences with unit increase looks the same as the
Catalan family tree, just with a different labelling and we obtain
an easy bijective correspondence between the self-describing
sequences and the sequences with unit increase. We could use this
bijective connection to show that the Catalan numbers count the
number of self-describing sequences. Instead, we provided a direct
proof of Theorem~\ref{thm:count} and the reason is that there is
an important difference in the distribution of labels in the
Catalan family tree and the tree of the sequences with unit
increase.

\begin{theorem}
For $r=0,\dots,n$, the number of $n$-th generation vertices in the
tree of sequences with unit increase labelled by $r$ is
\[ \frac{r+1}{n+1}\binom{2n-r}{n}. \]
\end{theorem}
\begin{proof}
Let $a=(a_0,a_1,\dots,a_n)$ be a sequence with unit increase.
Following Exercise 6.19.u in \cite{stanley:ec2}, we define, for
$i=0,\dots,n-1$,
\[ b_i = a_i - a_{i+1} + 1. \]
Construct a sequence of $n$ $1$'s and $n-a_n$ negative $1$'s by
replacing each $b_i$, $i=0,\dots,n-1$ by one $1$ followed by $b_i$
negative $1$'s. The newly obtained sequence has non-negative
partial sums. The correspondence between the sequences in $\A_n$
with unit increase that end by $r$ and the sequences of $n$ $1$'s
and $n-r$ negative $1$'s with non-negative partial sums is
bijective. It is shown in \cite{bailey:catalan} that the number of
sequences with non-negative partial sums that consist of $n$ $1$'s
and $k$ negative $1$'s is equal to
\[ \frac{n+1-k}{n+1}\binom{n+k}{n}\]
and this implies our claim.
\end{proof}

In passing, we make a slightly more general remark. Namely, for a
fixed positive integer $m$, consider the sequences with the
property that $a_0=0$ and $0 \leq a_{i+1} \leq a_i+ m$, for all
indices. Such sequences are called \emph{sequences with
$m$-increase}. We can easily construct the rooted labelled tree
that corresponds to such sequences. For a sequence
$(a_0,a_1,\dots,a_n)$ with $m$-increase, define, for
$i=0,\dots,n-1$,
\[ b_i = a_i - a_{i+1} + m. \]
Following the same approach as before, construct a sequence of $n$
$m$'s and $n-a_n$ negative $1$'s by replacing each $b_i$,
$i=0,\dots,n-1$ by one $m$ followed by $b_i$ negative $1$'s. The
newly obtained sequence has non-negative partial sums and the
correspondence between the sequences $(a_0,a_1,\dots,a_n)$ with
$m$-increase that end by $r$ and the sequences of $n$ $1$'s and
$mn-r$ negative $1$'s with non-negative partial sums is bijective.
Such sequences are discussed in \cite{frey-s:bailey}, where simple
recursive formulae for their number is provided. Unfortunately,
closed formulae are not provided yet, but we note that the number
of $n$-th generation sequences with $m$-increase is given by
$c_m(n+1)$ where
\[ c_m(n) = \frac{1}{mn+1}\binom{(m+1)n}{n}.  \]
The last displayed number is the generalization of the Catalan
numbers which counts, for example, the number of rooted
$(m+1)$-ary trees with $n$ interior vertices.

It is worth nothing that Julian West \cite{west:catalantree}
recursively constructs a rooted labelled tree whose root is
labelled by $2$ and each vertex labelled by $x$ has $x$ children
labelled by $2,3,\dots,x+1$. This tree, which West calls a Catalan
tree, looks again exactly like the Catalan family tree, but with
different labels. In fact, the tree of the sequences with unit
increase can be obtained from the Catalan tree constructed by
Julian West by decreasing all labels by 2.

Similarly, in the spirit of the Julian West construction, for any
positive integer $m$, construct a rooted labelled tree whose root
is labelled by $m+1$ and each vertex labelled by $x$ has $x$
children labelled by $m+1,m+2,\dots,m+x$. The tree of sequences
with $m$-increase can be obtained from this tree by decreasing all
labels by $m+1$.

\subsection*{Mirror symmetry and mutually describing sequences}
We introduce another endomorphism $\gamma: \A \to \A$ transforming
sequences in $\A$ by
\[ (\gamma a)_i = \# \{j\;|\; j < i,\; a_j \geq a_i \}. \]
Clearly $\gamma=\mu\delta$ where $\mu$ is the \emph{mirror
involution} of $\A$ given by
\[ (\mu a)_i = i-a_i. \]
We call $\mu$ the mirror involution of $\A$ since $\mu$ mirrors
the tree $\T$ through its vertical axis of symmetry.

The endomorphism $\gamma$ is studied in \cite{sunik:mds}. Clearly,
$\gamma$ has no fixed points other than the sequence $(0)$.
However, $\gamma$ has a lot of double points. If $a$ is a double
point of $\gamma$ then so is $b=\gamma a$. Moreover, then $\gamma
b= a$ and the sequences $a$ and $b$ mutually describe each other.

\begin{theorem}[\cite{sunik:mds}]
Repeated applications of $\gamma$ to any sequence in $\A$
eventually produce a double point of $\gamma$. The number of
application needed to reach a double point in $\A_n$ is $O(n^2)$
and there are more than $2^n$ such points.
\end{theorem}

The sequence that counts the number of double points of $\gamma$
in the $n$-th generation starts as follows
\[ 1, 2, 4, 10, 26, 70, 216, \dots \]
This sequence does not appear in the Encyclopedia of Integer
Sequences \cite{sloane-p:encyclopedia} nor in the online version
\cite{sloane:online} as of January 2002. It is interesting that
we have such a good understanding of the fixed points of
$\delta$, via the Catalan family tree, but we are still not able
to count the number of double points of the mirror related
endomorphism $\gamma=\mu\delta$.

Some other endomorphisms leading to fixed or double points are
studied in \cite{sunik:mds}. For one of them, the set of double
points of length $n$ is in bijective correspondence with the Young
tableaux of size $n$.

\subsection*{Acknowledgements}
Thanks to Richard Stanley and Louis Shapiro for their interest
and input.

\bibliographystyle{amsalpha}

\end{document}